\documentclass[12pt,twoside]{article}
\usepackage{graphicx}
\newtheorem{theorem}{\bf Theorem}[section]

\newtheorem{lemma}[theorem]{\bf Lemma}
\newtheorem{corollary}[theorem]{\bf Corollary}

{%{\footnotesize\sc }}
\input{amssym}

\date{September 10, 2009}
\begin{document}
\title{{\large A generalization of Ramsey theory for stars and one matching}}
\author{\small A. Khamseh$^{\textrm{a}}$, G.R. Omidi$^{\textrm{a},\textrm{b},1}$\\
\small  $^{\textrm{a}}$Department of Mathematical Sciences, Isfahan
University
of Technology,\\ \small Isfahan, 84156-83111, Iran\\
\small  $^{\textrm{b}}$School of Mathematics, Institute for Research
in Fundamental Sciences (IPM),\\
\small  P.O. Box: 19395-5746, Tehran,
Iran\\
\small \texttt{E-mails:khamseh@math.iut.ac.ir, romidi@cc.iut.ac.ir}}
\date {}

\maketitle\footnotetext[1] {\tt This research was in part supported
by a grant from IPM (No.89050037)} \vspace*{-0.5cm}

\begin{abstract}
\footnotesize
A  recent question in generalized Ramsey theory is that
for fixed positive integers $s\leq t$,
at least how many vertices can be covered by the vertices
of no more than $s$ monochromatic members of the family
$\cal F$ in every edge coloring of
$K_n$ with $t$ colors. This is related to
{{$d$-chromatic Ramsey numbers}} introduced by
Chung and Liu. In this paper, we first
compute these numbers for stars generalizing
the well-known result of  Burr and Roberts.
Then we extend  a result of
Cockayne and Lorimer to compute $d$-chromatic Ramsey numbers
for stars and one matching.
%%%%%%%%%%%%%%%%%%%%%%%%%%%%%
\\\\{\small {Keywords}: $d$-chromatic Ramsey number, edge coloring.}
\\\\{\small
{AMS subject classification}: 05C55, 05D10.}

\end{abstract}
%%%%%%%%%%%%%%%%%%%%%%%%%%%%%%%%%%%%%%%%%%%%%%%%%%%%%%%%%%%%%%%%%%%%%%%%
\section{\normalsize{Introduction}}
\small

%\section{{Introduction}}
Ramsey theory is an area
of combinatorics  which uses techniques from many branches of
mathematics and is currently among the most active areas in
combinatorics. Let $G_1, \ldots , G_c$ be graphs. The{\it{ Ramsey
number}} denoted by $r(G_1, \ldots , G_c)$ is defined to be the
least number $p$ such that if the edges of the complete graph
$K_p$ are arbitrarily colored with $c$ colors, then for some $i$
the spanning subgraph whose edges are colored with the $i$-th
color contains  $G_i$. More information about the Ramsey numbers
of known graphs can be found in the survey \cite{sur}.

\bigskip
There are
various types of Ramsey numbers that are important in the study of
classical Ramsey numbers and also hypergraph Ramsey numbers.
A question  recently proposed
by  Gy\'{a}rf\'{a}s et al. in \cite{gyarfas};
for fixed positive integers $s\leq t$,
at least how many vertices can be covered by the vertices
of no more than $s$ monochromatic members of the family
$\cal F$ in every edge coloring of
$K_n$ with $t$ colors.  Several problems and interesting
conjectures was presented in \cite{gyarfas}.
A basic problem here is to find
the largest $s$-colored element of $\cal F$
that can be found in every $t$-coloring of $K_n$.
The answer for  matchings when $s=t-1$ was given in \cite{gyarfas};
every $t$-coloring of $K_n$ contains a $(t-1)$-colored matching
of size $k$ provided that $n\geq 2k+[\frac{k-1}{2^{t-1}-1}]$.
Note that for $t=2,3,4$, we can guarantee the existence
of a $(t-1)$-colored path on ${2k}$ vertices instead of  a matching
of size $k$.
This was proved  in \cite{gerencser}, \cite{meenakshi}, and
\cite{khomidiforests}, respectively.

\bigskip
The above mentioned question is  related to an old problem
of Chung and Liu \cite{chli2};
for a given graph $G$ and for fixed  $s$, $t$, find the smallest $n$ such that in
every $t$-coloring of the edges of $K_n$ there is a copy of $G$ colored with at most $s$
colors.
More generally, let $1\leq d<c$ and let $t={c \choose d}$.
Assume that $A_1, \ldots, A_t$ are all $d$-subsets of a set
containing $c$ distinct colors. Let $G_1, \ldots, G_t$ be graphs.
The {\it{$d$-chromatic Ramsey numbers}} denoted by $r_d^c(G_1,  \ldots,
G_t)$ is the least number $p$ such that, if the edges
of the complete graph $K_p$ are arbitrarily colored  with $c$
colors, then for some $i$, the subgraph whose edges are colored by
colors in $A_i$
contains $G_ i$.

\bigskip
For complete graphs these numbers was partially
determined in \cite{chli2} and \cite{harborth}.
However for these graphs, the problem is very few known
and there are many open problems.
For stars, when  $d=1$ it is a well-know
result \cite{ramseyforstars}, and for $d=t-1=2$
the  value of $r_2^3(K_{1, i}, K_{1, j}, K_{1, l})$
was determined in   \cite{chli1}. For stars and one matching,
when $d=1$ it is again a well-known result; see \cite{cockayne}.

\bigskip
In this note, we first extend  the result of  \cite{chli1}
for stars to
arbitrary $c,d$ with $d=c-1\geq 2$. Then we replace  one of the
stars by a matching generalizing the result of
Cockayne and Lorimer to any  $c,d$ with $d=c-1\geq 2$.
To fix the notation, we use $r_{t-1}^{ t}(G_1,\ldots, G_t)$
to denote the minimum $p$ such that any
coloring of the edges of $K_p$ with $t$ colors $1,\ldots,t$
contains  a copy of $G_i$ for some $i$, missing the color $i$.
It is assumed throughout
the paper that $m_i\leq m_j$, where $i\leq j$ and
graphs are all simple and finite. A matching of size $m$
is denoted by $mP_2$ and a star of order $m+1$ by $K_{1,m}$.

%A colored graph $G$ is a complete graph with certain number of
%vertices whose edges are colored with four colors $1$, $2$, $3$,
%and $4$.

%%%%%%%%%%%%%%%%%%%%%%%%%%%%%%%%%%%%%%%%%%%%%%%%%%%%%%%%%%%%%%%%%%%%%%%%%

\section{{\normalsize{$(t-1)$-colored stars in $t$-colored complete graphs}}}\label{secstars}

In this section, we denote $\sum_{i=1}^{t}(m_i-1)$ briefly by
$\sum$. Let ${\rm{ex}}(p, H)$ be the maximum number of edges in a
graph on $p$ vertices which is $H$-free, i.e. it does not have $H$
as a subgraph. It is easily seen that ${\rm{ex}}(p, K_{1, m})\leq
\frac{p(m-1)}{2}$. We use this fact in the proof of Theorem
\ref{upperstar}.

\begin{theorem}\label{upperstar}
Let $x=\left[\frac{\sum+t-1}{t-1}\right]$. Then $r_{t-1}^{
t}(K_{1, m_1},\ldots,K_{1, m_t}) \leq x+1.$
\end{theorem}
{\it Proof.} Consider an edge coloring of $K_{x+1}$ with $t$
colors $1,\ldots,t$. Let $l_i$, $1\leq i\leq t$ be the number of
edges in color $i$ and $l=\sum_{i=1}^t l_i$. If for every $i$, we
have  $l-l_i\leq \frac{(x+1)(m_i-1)}{2}$, then  $x+1\leq
\frac{\sum+t-1}{t-1}$, a contradiction. So there exists an $i$
with $l-l_i>\frac{(x+1)(m_i-1)}{2}$. Hence the induced subgraph on
the edges with colors $\{1,\ldots,t\}-i$  contains a
$K_{1,m_i}$, as required.$\hfill \dashv $ \\

For graphs $G_1$, $G_2$, and $G_3$ with $|G_1|\leq |G_2|\leq
|G_3|$ it is shown \cite{chli2} that $r_{2}^{ 3}(G_1, G_2,
G_3)\leq r(G_1,G_2)$ and the equality holds if $|G_3| \geq
r(G_1,G_2)$, where $|G|$ is the number of vertices of $G$. Note
that for graphs $G_1$ and $G_2$, $r(G_1,G_2)=r_{1}^{2}(G_1,G_2)$.
So we can replace $|G_3| \geq r(G_1,G_2)$ by $|G_3| \geq
r_{1}^{2}(G_1,G_2)$. Theorem \ref{general}, is a trivial
generalization of this result.

%This result can be easily generalized the the following theorem.

\begin{theorem}\label{general}
Let $G_1, \ldots,G_t$ be graphs.
%with $|G_1|\leq \ldots\leq
%|G_t|$.
Then we have $r_{t-1}^{t}(G_1,\ldots,G_t)\leq
r_{t-2}^{t-1}(G_1,\ldots,G_{t-1})$ and the equality holds if
$|G_t|\geq r_{t-2}^{t-1}(G_1,\ldots,G_{t-1})$.
\end{theorem}
{\it Proof.} Let $l=r_{t-2}^{ t-1}(G_1,\ldots,G_{t-1})$ and
$c:E(G)\rightarrow \{1,2,\ldots,t\}$ be a coloring of $G=K_l$.
Define a new coloring $c'$ of $G$ with $t-1$ colors
${\texttt{1}},{\texttt{2}}, \ldots, {\texttt{t-1}}$ with
$c'(e)=\texttt{i}$ if $c(e)=i$, $1\leq i\leq t-2$,  and
$c'(e)=\texttt{t-1}$ if $c(e)=t-1$ or $c(e)=t$. By definition, $G$
contains a copy of $G_i$, for some $1\leq i\leq t-1$, in colors
$\{{\texttt{1}},\ldots,{\texttt{t-1}}\}-\{{\texttt{i}}\}$ which
implies that $G$ contains a copy of $G_i$, for some $1\leq i\leq
t$, in colors $\{1,\ldots,t\}-\{i\}$, as required.

\bigskip
Now suppose that $|G_t|\geq r_{t-2}^{t-1}(G_1,\ldots,G_{t-1})$. By
definition, there exists a coloring of $K_{l-1}$ with $t-1$ colors
such that $K_{l-1}$ does not contain $G_i$, for some $1\leq i\leq
t-1$, in colors $\{1,\ldots,t-1\}-\{i\}$. This is also a coloring
of $K_{l-1}$ with $t$ colors without $G_i$, $1\leq i\leq t$, in
colors $\{1,\ldots,t\}-\{i\}$. Thus
$$l-1<r_{t-1}^{t}(G_1,\ldots,G_t)\leq l=r_{t-2}^{t-1}(G_1,\ldots,G_{t-1}),$$
completing the proof.$\hfill \dashv $ \\

For abbreviation, we let $R_t=r_{t-1}^{ t}(K_{1, m_1},\ldots,K_{1,
m_t})$ and $x_t=\left[\frac{(\sum_{i=1}^tm_i)-1}{t-1}\right]$.
Then by Theorem \ref{general}, we can assume that $m_t+1\leq
R_{t-1}$. On the other hand, $R_t\leq R_{t-1}\leq \ldots\leq R_2$
and by Theorem \ref{upperstar}, $R_t\leq x_t+1$. Hence $m_t\leq
R_{t-1}-1\leq x_{t-1}$, which implies that $(t-2)m_t\leq
(\sum_{i=1}^{t-1}m_i)-1$. The last inequality is equivalent to
$x_{t-1}\geq x_t$. Similarly, $m_{t-1}\leq m_t\leq R_{t-1}-1\leq
R_{t-2}-1\leq x_{t-2}$ implies $x_{t-2}\geq x_{t-1}$. We continue
in this way, obtaining that $x_i\leq x_j$ for $j<i$. Using this
observation, we next find a lower bound for $r_{t-1}^{ t}(K_{1,
m_1},\ldots,K_{1, m_t})$.

\begin{theorem}\label{lowerstar}
Let $x=\left[\frac{\sum+t-1}{t-1}\right]$ and $m_t\leq R_{t-1}-1$.
Then
$$r_{t-1}^{ t}(K_{1, m_1},\ldots,K_{1, m_t})> x-1.$$
\end{theorem}
{\it Proof.} Let $p=x-\epsilon$ where $\epsilon=1$ if $x$ is odd
and $\epsilon=0$, otherwise. By Vizing's Theorem, there exists a
proper edge coloring of $K_p$ with $p-1$ colors. Let $r$, $1\leq
r< t$  be the smallest index such that $p-m_r\geq 0$ and
$p-m_{r+1}<0$ if it exists, and $r=t-1$ otherwise. Partition these $p-1$ colors
into $r+1$ new color classes as follows. Consider $p-m_i$ colors
as the new color $i$, for $1\leq i\leq r$ and all of the remaining
colors as the new color $r+1$. Note that since $p\leq x=x_t\leq
x_r=\left[\frac{(\sum_{i=1}^rm_i)-1}{r-1}\right]$, we have
$\sum_{i=1}^r(p-m_i)\leq p-1$.  This yields an edge coloring of
$K_p$ with $t$ colors $\{1,\ldots,t\}$ such that for each $i\leq
r$, every vertex $v$ is adjacent to at least $p-m_i$ edges in
color $i$ which  rules out the existence of $K_{1,m_i}$ in colors
$\{1,\ldots,t\}-\{i\}$. Moreover for $i\geq r+1$, no $K_{1,m_i}$
occurs since $p<m_{i}$. Hence $r_{t-1}^{ t}(K_{1,
m_1},\ldots,K_{1, m_t})> p$, which is
our assertion.$\hfill \dashv $ \\

The above proof  gives more, namely if
$x=\left[\frac{\sum+t-1}{t-1}\right]$ is even, then $$r_{t-1}^{
t}(K_{1, m_1},\ldots,K_{1, m_t})> x.$$ Combining this with Theorem
\ref{upperstar}, we conclude the following.

\begin{corollary}\label{xeven}
Let $x=\left[\frac{\sum+t-1}{t-1}\right]$ be even and $m_t\leq R_{t-1}-1$. Then
$$r_{t-1}^{ t}(K_{1, m_1},\ldots,K_{1, m_t})=x+1.$$
\end{corollary}

\vspace{.5cm} \noindent {\bf Remark}. Let  $v_1, \ldots, v_x$ be
vertices of $K_x$, where $x$ is odd. Eliminating $v_x$, there
exists corresponding matching $M_{v_x}$ %corresponding to $v_x$
containing $(x-1)/2$ parallel edges $v_1v_{x-1}, v_2v_{x-2}
,\ldots, v_{(x-1)/2}v_{(x+1)/2}$. Order these edges as above.
Similarly, for each vertex $v_i$, $1\leq i\leq x-1$, there exists
a matching $M_{v_i}$ containing $(x-1)/2$ ordered edges. These
matchings are used to construct certain edge colorings of $K_x$,
for example as in the proof of Theorem \ref{stars}.

\begin{theorem}\label{stars}
Let $x=\left[\frac{\sum+t-1}{t-1}\right]$, $m_t\leq R_{t-1}-1$ and $\sum=q(t-1)+h$,
where $0\leq h\leq t-2$. Then
$$r_{t-1}^{ t}(K_{1, m_1},\ldots,K_{1, m_t})=\left \{ \begin{array}{cc} x &~~~ {\rm if} ~x ~
{\rm{is~ odd}}, h=0 ~{\rm{and }} ~{\rm{some~}}  m_i {\rm{~is ~even }},  \\
 x+1 & {\rm otherwise.}  \end{array}\right. $$
\end{theorem}
{\it Proof.} If $x$ is even, then by Corollary \ref{xeven},
$r_{t-1}^{ t}(K_{1, m_1},\ldots,K_{1, m_t})=x+1$. So we may assume
that $x$ is odd.  We consider three cases as follows.

\begin{figure}[h]
  \begin{center}
  \includegraphics[width=5.5cm]{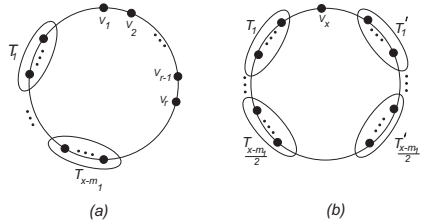}%\hspace*{2cm}
  \caption{Graph $K_x$ }\label{G1}
  \end{center}
\end{figure}

\vspace{.5cm} \noindent{\bf  {Case $1$.}} $h\geq 1$. Then
$r=x+\sum+t-tx=\sum+t-(t-1)x\geq 2$. Partition the vertices of
$K_x$ as $v_1, v_2, \ldots, v_r$ plus $x-m_1$ classes $T_1,
\ldots, T_{x-m_1}$ such that for $1\leq i\leq t$, we have
$T_i=\{u_{ij} : 1\leq j\leq n_i\}$, where $n_i$ is the largest
value $\lambda$ for which $i\leq x-m_\lambda$.
%$x-m_i$ first classes contains vertex labelled $u_{i}$. Then
%classes $T_1, \ldots, T_{x-m_t}$ contains vertices labelled
%$u_1,\ldots, u_t$, and classes $T_{x-m_t+1}, \ldots,
%T_{m_t-m_{t-1}}$ contains vertices labelled $u_1, \ldots,
%u_{t-1}$, etc.
For each vertex $u_{ij}$, $1\leq j\leq t$, paint with $j$ all
edges in $M_{u_{ij}}$. Let $v_1$ and $v_r$ be the vertices next to
$T_1$ and $T_{x-m_1}$, respectively (see Figure \ref{G1}($a$)).
For the vertex $v_1$ (respectively $v_r$) paint the edge
$e=u_{ij}v_l\in M_{v_1}$ (respectively $M_{v_r}$) with $j$ and
paint the edge $e=u_{ij}u_{i'j'}\in M_{v_1}$ (respectively
$M_{v_r}$) with $j$ if either $i<i'$ or $i=i'$ and $j<j'$
(respectively if either $i>i'$ or $i=i'$ and $j>j'$).
%, where $u_i\in T_j$, for $j\in \{1,\ldots,[\frac{x-m_1}{2}]+1\}$
%(also for $j\in \{[\frac{x-m_1}{2}]1,\ldots,x-m_1\}$).
The result is an edge coloring of $K_x$ with the property that for
each vertex, every color $i$ appears on at least $x-m_i$ edges;
that is, $r_{t-1}^{ t}(K_{1, m_1},\ldots,K_{1, m_t})>x$,  and so
by Theorem \ref{upperstar}, our assertion follows.

\vspace{.5cm} \noindent{\bf  Case $2$.} $h=0$, and every $m_i$ is
odd. Then $\sum=q(t-1)$, and $(t-1)(q-x)+t=1$. Partition the
vertices of $K_x$ as a single vertex $v_x$ plus $(x-m_1)/2$
classes $T_1, \ldots, T_{(x-m_1)/2}$,  and $(x-m_1)/2$ classes
$T'_1, \ldots, T'_{(x-m_1)/2}$ such that $T_i=\{u_{ij} : 1\leq
j\leq n_i\}$ and $T'_i=\{u'_{ij} : 1\leq j\leq n_i\}$, where $n_i$
is the largest value $\lambda$ for which $2i\leq x-m_\lambda$. Set
the classes $T_1, \ldots, T_{(x-m_1)/2}$ one side of $v_x$ and the
classes $T'_1, \ldots, T'_{(x-m_1)/2}$ on the other side of $v_x$,
respectively (see Figure \ref{G1}($b$)).
%$T_1, \ldots, T_{(x-m_t)/2}$,
%and also $T'_1, \ldots, T'_{(x-m_t)/2}$ contains $t$ vertices
%labelled $u_{1}, \ldots, u_{t}$, and $T_{(x-m_t)/2+1}, \ldots,
%T_{(m_t-m_{t-1})/2}$, and also $T'_{(x-m_t)/2+1}, \ldots,
%T'_{(m_t-m_{t-1})/2}$ contains $t-1$ vertices labelled $u_{1},
%\ldots, u_{t-1}$, etc. Set the classes $T_1, \ldots,
%T_{(x-m_1)/2}$ one side of $v_x$ and the classes $T'_1, \ldots,
%T'_{(x-m_1)/2}$ on the other side of $v_x$, respectively.
For each vertex $u_{ij}$ (also $u'_{ij}$), $1\leq j\leq t$, paint
with $j$ all edges in $M_{u_{ij}}$ (also $M_{u'_{ij}}$). Moreover,
for the vertex $v_x$, paint with $j$ the edge $e=u_{ij}u'_{ij}\in
M_{v_x}$. The result is an edge coloring of $K_x$ with the
property that for each vertex, every color $i$ appears on exactly
$x-m_i$ edges; that is, $r_{t-1}^{ t}(K_{1, m_1},\ldots,K_{1,
m_t})>x$,  and so by Theorem \ref{upperstar}, our assertion
follows.

\vspace{.5cm} \noindent{\bf  Case $3$.} $h=0$, and some $m_i$ is
even. Let $m_{i_0}$ be even. Then $x-m_{i_0}$ is odd. Suppose,
contrary to our claim, that  $r_{t-1}^{ t}(K_{1, m_1},\ldots,K_{1,
m_t})> x$. Consider the correspondent edge coloring of $K_x$ with
$t$ colors $1,\ldots,t$. As a sufficient condition, the degree of
each vertex in color $i$, $1\leq i\leq t$, is exactly $x-m_i$.
Then the induced subgraph with the edges in color $i_0$, is
$(x-m_{i_0})$-regular on $x$ vertices, a contradiction. Hence
$r_{t-1}^{ t}(K_{1, m_1},\ldots,K_{1, m_t})\leq x$, and so by
Theorem \ref{lowerstar}, our assertion follows.$\hfill \dashv $ \\

It may be worth reminding the reader that  Theorem \ref{stars} is
consistent with the well-known result of \cite{ramseyforstars}
that $r(K_{1, n},K_{1, m})=m+n-\epsilon$ where $\epsilon=1$  if
both $n$ and $m$ are even and $\epsilon=0$, otherwise.

%As a corollary, we have the value of standard Ramsey number
%$r_2(K_{1, n},K_{1, m})$ (see \cite{ramseyforstars}).

%\begin{corollary}
%$r_2(K_{1, n},K_{1, m})=m+n-\epsilon$  where $\epsilon=1$  if both
%$n$ and $m$ are even and $\epsilon=0$, otherwise.
%\end{corollary}
%%%%%%%%%%%%%%%%%%%%%%%%%%%%%%%%%%%%%%%%%%%%%%%%%%%%%%%%%%%%%%%%%%%%%%%%%%%%%%%%%
%%%%%%%%%%%%%%%%%%%%%%%%%%%%%%%%%%%%%%%%%%%%%%%%%%%%%%%%%%%%%%%%%%%%%%%%%%%%%%%%%%%%%%%%%%
\section{{\normalsize{$(t-1)$-colored stars-matching in $t$-colored complete graphs}}}\label{secstarstripes}

%This result can be easily generalized the the following theorem.

 In this section, we calculate
$r_{t-1}^{t}(K_{1,m_1},\ldots,K_{1,m_{t-1}},sP_2)$. In
\cite{cockayne} the value of
$r_1^2(K_{1,m_1},sP_2)=r(K_{1,m_1},sP_2)$ has been determined, so
we can assume that $t\geq 3$. For abbreviation, we write $R$
instead of $r_{t-1}^{t}(K_{1,m_1},\ldots,K_{1,m_{t-1}},sP_2)$ and
denote $\sum_{i=1}^{t-1}(m_i-1)$ briefly by $\sum$.
If $2s\geq r_{t-2}^{t-1}(K_{1,m_1},\ldots,K_{1,m_{t-1}})$, %or
%$m_{t-1}+1\geq R_{t-2,t-1}(K_{1,m_1},\ldots,K_{1,m_{t-2}},sP_2)$
then by Theorem \ref{general},
$R=r_{t-2}^{t-1}(K_{1,m_1},\ldots,K_{1,m_{t-1}})$. %, or
%$R=R_{t-2,t-1}(K_{1,m_1},\ldots,K_{1,m_{t-2}},sP_2)$,
%respectively.
Therefore in the following two lemmas we assume
$2s<r_{t-2}^{t-1}(K_{1,m_1},\ldots,K_{1,m_{t-1}})$. %and
%$m_{t-1}+1<R_{t-2,t-1}(K_{1,m_1},\ldots,K_{1,m_{t-2}},sP_2)$.

\begin{lemma}\label{R=2s}
If $t\geq 3$, $\sum<(2t-3)s-t+2$, and
$2s<r_{t-2}^{t-1}(K_{1,m_1},\ldots,K_{1,m_{t-1}})$, then $R=2s$.
\end{lemma}
{\it Proof.} Since
$2s<r_{t-2}^{t-1}(K_{1,m_1},\ldots,K_{1,m_{t-1}})$, there exists
an edge coloring of $K_{2s-1}$ with colors $1,\ldots, t-1$, such
that for each $i$, $1\leq i\leq t-1$, the induced subgraph on the
edges with colors $\{1,\ldots, t-1\}-i$ does not contain
$K_{1,m_i}$. This also can be considered as an edge coloring of
$K_{2s-1}$ with $t$ colors $1,\ldots, t$ such that in addition,
the induced subgraph on the edges with colors $\{1,\ldots, t-1\}$
does not contain $sP_2$; that is, $R>2s-1$.

\bigskip
We now show that $R\leq 2s$. Consider an edge coloring of $K_{2s}$
with colors $1,\ldots,t$. Let $M$ be the maximal matching of edges
with colors $1,\ldots,t-1$. Then $M$ has at most $s'\leq s-1$
independent edges, since otherwise we are done. Let $W$ be the set
of those vertices that are not incident with these $s'$ edges.
Note that $|W|\geq 2$, and every edge incident with two vertices
in $W$ has color $t$. Moreover, every vertex is incident with at
least $2s-m_i$ edges in color $i$, $1\leq i\leq t-1$, since otherwise we are done. Thus every
vertex is incident with at least $2(t-1)s-\sum-(t-1)$ edges in
colors $1,\ldots,t-1$. Since $\sum<(2t-3)s-t+2$, each of the
vertices $w_1,w_2\in W$ is incident with at least $s$ edges in
colors $1,\ldots,t-1$; that is, there exists $e=uv\in M$ such that
the color of both $w_1u$, and $w_2v$ belong to $
\{1,\ldots,t-1\}$,
which contradicts the maximality of $M$.$\hfill \dashv $ \\

\begin{lemma}\label{R=sigma+s/t-1}
If $t\geq 3$, $\sum\geq (2t-3)s-t+2$, and
$2s<r_{t-2}^{t-1}(K_{1,m_1},\ldots,K_{1,m_{t-1}})$, then
$R=\left\lceil\frac{\sum+s}{t-1}\right\rceil+1$.
\end{lemma}
{\it Proof.} Let $l=\left\lceil\frac{\sum+s}{t-1}\right\rceil$. To
prove $R\leq l+1$, consider an edge coloring of $K_{l+1}$ with $t$
colors $1,\ldots,t$. Let $M$ be the maximal matching of edges with
colors $1,\ldots,t-1$. Then $M$ has at most $s'\leq s-1$
independent edges, since otherwise we are done. Let $W$ be the set
of those vertices that are not incident with these $s'$ edges.
Note that $|W|\geq 2$, and every edge incident with two vertices
in $W$ has color $t$. Moreover, every vertex is incident with at
least $l+1-m_i$ edges in color $i$, $1\leq i\leq t-1$. Thus every
vertex is incident with at least $(t-1)(l+1)-\sum-(t-1)$ edges in
colors $1,\ldots,t-1$. Let $w_1,w_2\in W$. Since
$l>\frac{\sum+s-1}{t-1}$, $(t-1)(l+1)-\sum-(t-1)>s-1$ and so each
of the vertices $w_1,w_2$ is incident with at least $s$ edges in
colors $1,\ldots,t-1$. Therefore, there exists $e=uv\in M$ such
that the color of both $w_1u$, and $w_2v$ belong to
$\{1,\ldots,t-1\}$, which contradicts the maximality of $M$.

\bigskip
We now turn our attention to the lower bound. Set $n_i=l-m_i$,
$1\leq i\leq t-1$. Partition the vertices of $K_l$ into $t-1$
classes $X_i$, $1\leq i\leq t-1$, with $|X_i|=n_i$ plus the set
$X$ consist of the rest of the vertices. Note that $n_i\geq 0$ and
$\sum_{i=1}^{t=1}n_i<l$. First let $z=\sum_{i=1}^{t-1}n_i$ be odd
and suppose that $x\in X$. By Vizing's Theorem, there exists an
edge coloring of the complete graph on $z+1$ vertices $\{x\}\cup
\bigcup_{i=1}^{t-1} X_i$ with $z$ colors. Set these $z$ colors
into $t-1$ color classes by considering $n_i$ colors as the new
color $i$, $1\leq i\leq t-1$. This yields an edge coloring of
$K_z$ with $t-1$ colors $\{1,\ldots,t-1\}$ such that every vertex
$v\in \{x\}\cup \bigcup_{i=1}^{t-1} X_i$ is adjacent to
$n_i=l-m_i$ edges in color $i$, $1\leq i\leq t-1$. Moreover, for
$1\leq i\leq t-1$, paint with $i$ the edges having  one vertex in
$X_i$ and one vertex in $X-\{x\}$. Finally, paint with $t$ all the
remaining edges. In this coloring of $K_l$, every vertex is
adjacent to at least $n_i$ edges in color $i$, $1\leq i\leq t-1$,
which rules out the existence of $K_{1,m_i}$ in colors
$\{1,\ldots,t\}-\{i\}$. Moreover, the subgraph on the edges with
colors $1,\ldots,t-1$ contains at most $s-1$ independent edges. We
now suppose that $z=\sum_{i=1}^{t-1}n_i$ is even. Let $x,y\in X$.
By Vizing's Theorem, there exists an edge coloring of the complete
graph on $z+2$ vertices $\{x,y\}\cup \bigcup_{i=1}^{t-1} X_i$ with
$z+1$ colors. Without loss of generality we can assume that $xy$
has color $1$. Partition these $z+1$ colors into $t-1$ color
classes by considering $n_1+1$ colors as the new color $1$ and
$n_i$ colors as the new color $i$, $2\leq i\leq t-1$. This yields
an edge coloring of $K_{z+2}$ with $t-1$ colors $\{1,\ldots,t-1\}$
such that every vertex $v\in \{x,y\}\cup \bigcup_{i=1}^{t-1} X_i$
is adjacent to at least $n_i=l-m_i$ edges in color $i$, $1\leq
i\leq t-1$. Moreover, for $1\leq i\leq t-1$, paint with $i$ the
edges having  one vertex in $X_i$ and one vertex in $X-\{x,y\}$.
Finally, paint with $t$ all the remaining edges and change the
color of $xy$ into $t$.
%Similar argument, with $x\in X$ replaced by $x,y\in X$, is
%applicable to obtain the coloring when $z=\sum_{i=1}^{t-1}n_i$ is
%even.
Again in this coloring of $K_l$, every vertex is adjacent to at
least $n_i$ edges in color $i$, $1\leq i\leq t-1$, which rules out
the existence of $K_{1,m_i}$ in colors $\{1,\ldots,t\}-\{i\}$.
Moreover, the subgraph on the edges with colors $1,\ldots,t-1$
contains at most $s-1$ independent edges.
Therefore, $R>l$, completing the proof.$\hfill \dashv $ \\

Combining Lemmas \ref{R=2s}, and \ref{R=sigma+s/t-1} with the
above discussion we have the following theorem.

\begin{theorem}
Let $t\geq 3$. Then
\begin{enumerate}
\item[$i.$] If $2s\geq R_{t-1}$, then $R=R_{t-1}$. \item[$ii.$] If
$2s< R_{t-1}$ and $\sum<(2t-3)s-t+2$, then $R=2s$. \item[$iii.$]
If $2s< R_{t-1}$ and $\sum\geq (2t-3)s-t+2$, then
$R=\left\lceil\frac{\sum+s}{t-1}\right\rceil+1$.
\end{enumerate}
\end{theorem}

\vspace{1.5cm}

 \end{document}